\documentclass[12pt]{amsart}
\usepackage{amssymb}
\usepackage{isolatin1}    
\newtheorem{thm}{Theorem}[section]
\newtheorem{lem}[thm]{Lemma}
\newtheorem{cor}[thm]{Corollary}
\newtheorem{prop}[thm]{Proposition}
\newtheorem{ex}[thm]{Example}

\newtheorem*{prob*}{Open problem}

\theoremstyle{definition}

\newtheorem{defi}[thm]{Definition}

\theoremstyle{remark}

\newtheorem{rem}[thm]{Remark}
\newtheorem*{rem*}{Remark}

\DeclareMathOperator{\s}{span}

\newcommand{\kringel}{\mathbin{\raise1pt\hbox{$\scriptstyle\circ$}}} 
\newcommand{\pkt}{\mathbin{\raise0pt\hbox{$\scriptstyle\bullet$}}}

\newcommand{\N}{\mathbb{N}}

\newcommand{\tr}{\mathop{\rm tr}}

\newcommand{\End}{\mathop{\rm End}}

\newcommand{\diag}{\mathop{\rm diag}}

\newcommand{\La}{\mathfrak{a}}
\newcommand{\Lb}{\mathfrak{b}}

\newcommand{\Lg}{\mathfrak{g}}

\newcommand{\CB}{\mathcal{B}}

\newcommand{\CH}{\mathcal{H}}

\newcommand{\im}{\mathop{\rm im}}

\newcommand{\al}{\alpha}
\newcommand{\be}{\beta}
\newcommand{\ga}{\gamma}
\newcommand{\de}{\delta}

\newcommand{\la}{\lambda}
\newcommand{\om}{\omega}

\newcommand{\Om}{\Omega}

\newcommand{\ra}{\rightarrow}  

\renewcommand{\phi}{\varphi}

\begin{document}

\title[$XA-AX=X^p$]{On the matrix equation {${\bf XA-AX=X^{p}}$}} 

\author[D. Burde]{Dietrich Burde}
\address{Fakult\"at f\"ur Mathematik\\
Universit\"at Wien\\
  Nordbergstr. 15\\
  1090 Wien \\
  Austria}
\date{\today}
\email{dietrich.burde@univie.ac.at}

\thanks{I thank Joachim Mahnkopf for helpful remarks.}

\begin{abstract}
We study the matrix equation $XA-AX=X^p$ in $M_n(K)$ for $1< p <n$. 
It is shown that every matrix solution $X$ is nilpotent and that the 
generalized eigenspaces of $A$ are $X$-invariant. For $A$ being a full 
Jordan block we describe how to compute all matrix solutions.
Combinatorial formulas for $A^mX^{\ell},X^{\ell}A^m$ and $(AX)^{\ell}$
are given.
The case $p=2$ is a special case of the algebraic Riccati equation.
\end{abstract}

\maketitle

\section{Introduction} 

Let $p$ be a positive integer. The matrix equation  
\begin{align*}
XA-AX & =X^p 
\end{align*}
arises from questions in Lie theory. In particular, the quadratic
matrix equation $XA-AX =X^2$ plays a role in 
the study of affine structures on solvable Lie algebras.\\
An {\it affine structure} on a Lie algebra $\Lg$ over a field $K$ is
a $K$-bilinear product $\Lg\times \Lg \ra \Lg$, $(x,y)\mapsto x\cdot y$
such that
\begin{align*}
x\cdot (y\cdot z)-(x\cdot y)\cdot z & = y\cdot (x\cdot z)-(y\cdot x)\cdot z\\
[x,y] & =x\cdot y -y\cdot x
\end{align*}
for all $x,y,z\in \Lg$ where $[x,y]$ denotes the Lie bracket of $\Lg$.
Affine structures on Lie algebras correspond to left-invariant affine
structures on Lie groups. They are important for affine manifolds and
for affine crystallographic groups, see  \cite{BU1}, \cite{GRS}, \cite{MIL}.\\
We want to explain how the quadratic matrix equations $XA-AX=X^2$
arise from affine structures.
Let $\Lg$ be a two-step solvable Lie algebra.
This means we have an exact sequence of Lie algebras
\begin{equation*}
0 \rightarrow \La \xrightarrow{\iota} \Lg \xrightarrow{\pi}
\Lb \rightarrow 0
\end{equation*}
with the following data: $\La$ and $\Lb$ are abelian Lie algebras,
$\phi : \Lb \mapsto \End (\La)$ is a Lie algebra representation,
$\Om \in Z^2(\Lb,\La)$ is a $2$-cocycle, and the Lie bracket of
$\Lg=\La \times \Lb$ is given by
 \begin{equation*}
[(a,x),(b,y)]:=(\phi (x)b-\phi(y)a+\Om(x,y),0).
\end{equation*}
Let $\om :\Lb \times \Lb \ra \La$ be a bilinear map and
$\phi_1,\, \phi_2 : \Lb \mapsto \End (\La)$ Lie algebra representations.
A natural choice for a left-symmetric product on $\Lg$ is the
bilinear product given by
\begin{equation*}
(a,x)\kringel (b,y):=(\phi_1(y)a+\phi_2 (x)b+\om(x,y),0).
\end{equation*}
One of the necessary conditions for the product to be left-symmetric
is the following:
\begin{align*}
\phi_1(x)\phi(y)-\phi(y)\phi_1(x) & = \phi_1(y)\phi_1(x).
\end{align*}

Let $(e_1,\ldots e_m)$ be a basis of $\Lb$ and write $X_i:=\phi_1(e_i)$
and $A_j:=\phi(e_j)$ for the linear operators. We obtain the matrix
equations

\begin{align*}
X_i A_j -A_j X_i & = X_jX_i
\end{align*}
for all $1\le i,j\le m$. In particular we have matrix equations of
the type $XA-AX=X^2$.

\section{General results}

Let $K$ be an algebraically closed field of characteristic zero.
In general it is quite difficult to determine the matrix solutions
of a nonlinear matrix equation. Even the existence of solutions is
a serious issue as illustrated by the quadratic matrix equation 
$$
X^2=\begin{pmatrix}
0 & 1 \\
0 & 0 
\end{pmatrix}
$$
which has no solution. On the other hand our equation $XA-AX=X^p$ always has
a solution, for any given $A$, namely $X=0$. However, if $A$ has a multiple eigenvalue,
then we have a lot of nontrivial solutions and there is no easy way to describe
the solution set algebraically. A special set of solutions is obtained
by the matrices $X$ satisfying $XA-AX=0=X^p$. First one can determine
the matrices $X$ commuting with $A$ and then pick out those satisfying
$X^p=0$.\\
Let $E$ denote the $n\times n$ identity.
We will assume most of time that $p\ge 2$ since
for $p=1$ we obtain the linear matrix equation $AX +X(E-A)=0$ which is
a special case of the Sylvester matrix equation $AX+XB=C$. 
Let $S: M_n(K)\ra M_n(K)$ with $S(X)=AX+XB$ be the Sylvester operator. 
It is well known that the linear operator $S$ 
is singular if and only if $A$ and $-B$ have a common eigenvalue, see
\cite{LAR}.
For $B=E-A$ we obtain the following result.

\begin{prop}
The matrix equation $XA-AX=X$ has a nonzero solution if and only if
$A$ and $A-E$ have a common eigenvalue.
\end{prop}

The general solution of the matrix equation
$AX=XB$ is given in \cite{GAT}.
We have the following results on the solutions of our general equation.

\begin{prop}
Let $A\in M_n(K)$. Then every matrix solution $X\in M_n(K)$ 
of $XA-AX=X^p$ is nilpotent and hence satisfies $X^n=0$.
\end{prop}

\begin{proof}
We have $X^k(XA-AX)=X^{k+p}$ for all $k\ge 0$. Taking the trace on both sides
we obtain $\tr (X^{k+p})=0$ for all $k\ge 0$. Let $\la_1,\ldots , \la_r$ be the
pairwise distinct eigenvalues of $X$. For $s\ge 1$ we have
$$
\tr (X^s)=\sum_{i=1}^r m_i\la_i^s.
$$
For $s\ge p$ we have $\tr (X^s)=0$ and hence
$$
\sum_{i=1}^r (m_i\la_i^p)\la_i^k = 0
$$
for all $k\ge 0$. This is a system of linear equations in the unknowns
$x_i=m_i\la_i^p$ for $i=1,2,\ldots ,r$. The determinant of its coefficients
is a Vandermonde determinant. It is nonzero since the $\la_i$ are pairwise
distinct. Hence it follows $m_i\la_i^p=0$ for all $i=1,2,\ldots ,r$. This 
means $\la_1=\la_2=\cdots =\la_r=0$ so that $X$ is nilpotent with $X^n=0$.
\end{proof}

Since for $p=n$ our equation reduces to $X^n=0$ and the linear matrix equation
$XA=AX$, we may assume that $p<n$.

\begin{prop}
Let $K$ be an algebraically closed field and $p$ be a positive integer.
If $X,A\in M_n(K)$ satisfy $XA-AX=X^p$ then $X$ and $A$ can be
simultaneously triangularized. 
\end{prop}

\begin{proof}
Let $V$ be the vector space generated by $A$ and all $X^{i}$. Since
$X$ is nilpotent we can choose a minimal $m\in \N$ such that $X^m=0$.
Then $V=\s \{ A,X,X^2,\ldots, X^{m-1}\}$.
We define a Lie bracket on $V$ by taking commutators. 
Using induction on $\ell$ we see that for all $\ell\ge 1$
\begin{equation}\label{a2}
X^{\ell}A-AX^{\ell}=\ell X^{p+\ell-1}
\end{equation}
Hence the Lie brackets are defined by 
\begin{align*}
[A,A] & = 0\\
[A,X^i] & = AX^i-X^iA=-iX^{p+i-1}\\
[X^i,X^j] & = 0.
\end{align*}
It follows that $V$ is a finite-dimensional 
Lie algebra. The commutator Lie algebra $[V,V]$ is abelian and
$V/[V,V]$ is $1$-dimensional. Hence $V$ is solvable.  
By Lie's theorem $V$ is triangularizable. Hence there is a basis such that 
$X$ and $A$ are simultaneously upper triangular.
\end{proof}

\begin{cor}
Let $X,A\in M_n(K)$ satisfy the matrix equation $XA-AX=X^p$. Then 
$A^iX^k A^j$ is nilpotent for all $k\ge 1$ and $i,j\ge 0$. 
So are linear combinations of such matrices.
\end{cor}

\begin{proof}
We may assume that $X$ and $A$ are simultaneously upper triangular.
Since $X$ is nilpotent, $X^k$ is strictly upper
triangular. The product of such a matrix with an upper
triangular matrix $A^i$ or $A^j$ is again strictly upper triangular.  
Moreover a linear combination of strictly upper triangular matrices is
again strictly upper triangular.
\end{proof}

\begin{prop}
Let $p\ge 2$ and $A\in M_n(K)$. If $A$ has no multiple eigenvalue
then $X=0$ is the only matrix 
solution of $XA-AX=X^p$. Conversely if $A$ has a multiple eigenvalue
then there exists a nontrivial solution $X\ne 0$.
\end{prop}

\begin{proof}
Assume first that $A$ has no multiple eigenvalue. Let $\CB=(e_1,\ldots , e_n)$
be a basis of $K^n$ such that $A=(a_{ij})$ and $X=(x_{ij})$ are upper
triangular relative to $\CB$. In particular $a_{ij}=0$ for $i>j$ and
$x_{ij}=0$ for $i\ge j$. Since all eigenvalues of $A$ are distinct, $A$ is
diagonalizable. We can diagonalize $A$ by a base change of the form
$e_i\mapsto \mu_1e_1+\mu_2e_2+\cdots + \mu_i e_i$ which also keeps
$X$ strictly upper triangular. Hence we may assume that $A$ is diagonal
and $X$ is strictly upper triangular. Then the coefficients
of the matrix $XA-AX=(c_{ij})$ satisfy 
$$
c_{ij}=x_{ij}(a_{jj}-a_{ii}),\quad x_{ij}=0 \; \text{for} \; i\ge j.
$$
Consider the lowest nonzero line parallel to the main diagonal in $X$.
Since $\al_{jj}-\al_{ii}\ne 0$ for all $i\ne j$ this line stays also nonzero
in $XA-AX$, but not in $X^p$ because of $p\ge 2$. It follows that $X=0$.\par
Now assume that $A$ has a multiple eigenvalue. There exists a basis of $K^n$
such that $A$ has canonical Jordan block form. Each Jordan block is an
matrix of the form
$$
J(r,\la)={\Small
\begin{pmatrix}
 \la & 1 & \hdots & 0 & 0\\
 0 & \la & \hdots & 0 & 0\\
\vdots & \vdots & \ddots & \vdots & \vdots\\
0 & 0  & \hdots & \la & 1\\
0 & 0  & \hdots & 0 & \la
\end{pmatrix}}\in M_r(K).
$$
For $\la=0$ we put $J(r)=J(r,0)$. It is $J(r,\la)=J(r)+\la E$.
Consider the matrix equation $XJ(r,\la)-J(r,\la)X=X^p$ in $M_r(K)$. 
It is equivalent to the equation
$XJ(r)-J(r)X=X^p$. If $r \ge 2$
it has a nonzero solution, namely the $r\times r$ matrix
$$
X={\Small
\begin{pmatrix}
0 &  \hdots & 1\\
\vdots & \ddots & \vdots \\
0  & \hdots & 0 
\end{pmatrix}}.
$$
Indeed, $XJ(r)-J(r)X=0=X^p$ in that case. 
Since $A$ has a multiple eigenvalue, it has a Jordan block of size
$r\ge 2$. After permutation we may assume that this is the first Jordan block
of $A$. Let $X\in M_r(K)$ be the above matrix and extend it to an $n\times n$-matrix
by forming a block matrix with $X$ and the zero matrix in $M_{n-r}(K)$.
This will be a  nontrivial solution of $XA-AX=X^p$ in $M_n(K)$. 
\end{proof}

\begin{lem}
Let $A,X\in M_n(K)$ and $A_1=SAS^{-1}$, $X_1=SXS^{-1}$ for some $S\in GL_n(K)$. 
Then $XA-AX=X^p$ if and only if $X_1A_1-A_1X_1=X_1^p$. 
\end{lem}

\begin{proof}
The equation $X^p=XA-AX$ is equivalent to 
\begin{align*}
X_1^p & =(SXS^{-1})^p = SX^pS^{-1} \\
& = (SXS^{-1})(SAS^{-1})-(SAS^{-1})(SXS^{-1})\\
& = X_1A_1-A_1X_1.
\end{align*}
\end{proof}
The lemma says that we may choose a basis of $K^n$ such that 
$A$ has canonical Jordan form.\\
Denote by $C(A)=\{S\in M_n(K) \mid SA=AS\}$ the centralizer of $A\in M_n(K)$.
Applying the lemma with $A_1=SAS^{-1}=A$, where $S\in C(A)\cap GL_n(K)$, we
obtain the following corollary.

\begin{cor}\label{cor7}
If $X_0$ is a matrix solution of $XA-AX=X^p$ then so is $X=SX_0S^{-1}$ for
any $S\in C(A)\cap GL_n(K)$.
\end{cor}

Let $\CB=(e_1,\ldots ,e_n)$ be a basis of $K^n$ such that
$A$ has canonical Jordan form. Then $A$ is a block matrix 
$$A=\diag (A_1,A_2,\ldots , A_k)$$
with $A_i\in M_{r_i}(K)$
and $A$ leaves invariant the corresponding subspaces of $K^n$.
Let $X$ satisfy $XA-AX=X^p$. 
Does it follow that $X$ is also a block matrix $X=\diag (X_1,\ldots, X_r)$ 
with $X_i\in M_{r_i}(K)$ relative to the basis $\CB$ ?
In general this is not the case.

\begin{ex}
The matrices
$$A=\begin{pmatrix} 0 & 0 & 0\\ 0 & 0 & 1\\ 0 & 0 & 0 \end{pmatrix},\quad
X=\begin{pmatrix} -1 & 0 & 1\\ -1 & 0 & 0\\ -1 & 0 & 1 \end{pmatrix}
$$
satisfy $XA-AX=X^2$.
\end{ex}
Here $A=\diag (J(1),J(2))$ leaves invariant the subspaces $\s \{e_1\}$
and $\s \{e_2,e_3\}$ corresponding to the Jordan blocks $J(1)$ and $J(2)$,
but $X$ does not. Also the subspace $\ker A$ is not
$X$-invariant. This shows that
the eigenspaces $E_{\la}=\{x\in K^n \mid Ax=\la x\}$ of $A$ 
need not be $X$-invariant.
However, we have the following result concerning the generalized
eigenspaces $\CH_{\la}=\{x\in K^n \mid (A-\la E)^kx=0 \;\text{for some} 
\;k\ge 0\}$ of $A$.

\begin{prop}
Let $A,X\in M_n(K)$ satisfy $XA-AX=X^p$ for $1<p<n$. Then the generalized
eigenspaces $\CH_{\la}$ of $A$ are $X$-invariant, i.e., $X \CH_{\la}\subseteq \CH_{\la}$.
\end{prop}

\begin{proof}
Let $\la$ be an eigenvalue of $A$ and  $\CH_{\la}$ be the generalized eigenspace.
We may assume that $A$ has canonical Jordan form such that $A=\diag (A_1,A_2)$
with $A_1=J(r,\la)$. We may also assume that $\la=0$. This follows by
considering $B=A-\la E$ instead of $A$ which satisfies $XB-BX=XA-AX=X^p$.
Let $v\in \CH_{0}$. Then there exists an integer $m\ge 0$ such that $A^mv=0$.
Let $r$ be an integer with $r\ge n$. We have $X^r=0$.
By induction on $k\ge 1$ we will show that
$$
A^{m+k-1}X^{r-k(p-1)}v=0 \quad \text{for} \quad 1\le k < \frac{r}{p-1}.
$$
This implies the desired result as follows: set
$r=1+ k(p-1)$. We can choose $k\ge 1$ such that 
$r\ge n$. Then $A^{m+k-1}Xv=0$ and hence $Xv\in \CH_0$.\\
For $k=1$ we have to show $A^mX^{r-(p-1)}v=0$.
By \eqref{a2} we have 
$$AX^{r-(p-1)}-X^{r-(p-1)}A=(p-1-r)X^r=0.$$ 
Hence $A$ and $X^{r-(p-1)}$ commute. It follows that also $A^m$ and $X^{r-(p-1)}$
commute. Hence $A^mX^{r-(p-1)}v=X^{r-(p-1)}A^mv=0$.\\
Assume now that $A^{m+k-2}X^{r-(k-1)(p-1)}v=0$. Then we have
$$
AX^{r-k(p-1)}-X^{r-k(p-1)}A = (k(p-1)-r)X^{r-(k-1)(p-1)}.
$$
Now we will use the following formula: let $s,\ell\ge 1$ be integers
and $A,X\in M_n(K)$ satisfying $XA-AX=X^p$, where $1<p<n$. Then there exist
integers $b_j=b_j(p,\ell,s)$ such that
\[
A^s X^{\ell} =\sum_{j=0}^s b_j X^{\ell+j(p-1)}A^{s-j}.
\]
This formula can be easily proved by induction. We will compute explicitly
the coefficients $b_j$ in the last section, see formula \eqref{a11}.
If we use the formula for $\ell=r-k(p-1)$ and $s=m+k-2$ then we obtain
\[
A^{m+k-2}X^{r-k(p-1)} =\sum_{j=0}^{m+k-2} b_j X^{r+(j-k)(p-1)}A^{m+k-2-j}.
\]
It follows 
\begin{align*}
A^{m+k-1}X^{r-k(p-1)}v & = A^{m+k-2}X^{r-k(p-1)}Av \\
 & = \sum_{j=0}^{m+k-2}b_j X^{r+(j-k)(p-1)}A^{m+k-j-1}v.
\end{align*}
Here all terms with $j\ge k$ vanish since $X^r=0$. On the other hand,
$A^{m+k-j-1}v=0$ for all $j\le k-1$. It follows $A^{m+k-1}X^{r-k(p-1)}v=0$.
\end{proof}

\begin{cor}
Let $A=\diag (A_1,A_2)$ be a block matrix such that $A_1\in M_r(K)$ and
$A_2\in M_s(K)$ have no common eigenvalues. Then any matrix solution
$X\in M_n(K)$ of $XA-AX=X^p$ for $1<p<n$ is a block matrix
$X=\diag (X_1,X_2)$ with $X_1\in M_r(K)$ and $X_2\in M_s(K)$ such that
$X_1A_1-A_1X_1=X_1^p$ and $X_2A_2-A_2X_2=X_2^p$.
\end{cor}

This says that looking at the solutions of $XA-AX=X^p$
we may restrict to the case that $A$ has exactly one eigenvalue 
$\la\in K$. Without loss of generality we may assume that $\la=0$.
We can say more on the solution set if $A$ has some particular
properties. 
The most convenient special case is that $A=J(n)$ is a full Jordan block.
Then we can determine all matrix solutions of $XJ(n)-J(n)X=X^p$.
This is done in the following section.

\section{The case $A=J(n)$}

We have already seen that $\ker A$ in general is not $X$-invariant.
However, it is true if $\ker A$ is $1$-dimensional. But this is the
case for $A=J(n)$.

\begin{lem}
Let $A,X\in M_n(K)$ satisfy $XA-AX=X^p$ for $1<p<n$ and assume that 
$\la$ is an eigenvalue of $A$ with $1$-dimensional
eigenspace $E_{\la}$ generated by $v\in K^n$. Then $Xv=0$. 
\end{lem}

\begin{proof}
We will show that $X^{\ell}v=0$ implies $X^{\ell-(p-1)}v=0$ for all
$\ell \ge p$. By \eqref{a2} we have
\begin{equation}
AX^{\ell -(p-1)}-X^{\ell -(p-1)}A= (p-1-\ell)X^{\ell}.
\end{equation}
Using $Av=\la v$ and $X^{\ell}v=0$ we obtain $(A-\la E)X^{\ell -(p-1)}v=0$
so that $X^{\ell -(p-1)}v \in E_{\la}=\s \{v\}$. Hence $X^{\ell -(p-1)}v=\mu v$ for
some $\mu \in K$. But since $X$ is nilpotent we have $\mu=0$.\par
Now we repeat this argument starting with $X^n v=0$. If we arrive
at $X^kv=0$ and $k\le p$ then $X^pv=0$ and in the next step $Xv=0$.
\end{proof}

\begin{prop}\label{prop12}
Let $A=J(n)$ and $X\in M_n(K)$ be a matrix solution of $XA-AX=X^p$ for
$1<p <n$. Then $X$ is strictly upper triangular.
\end{prop}

\begin{proof}
Let $(e_1,\ldots e_n)$ be the canonical basis of $K^n$. Then $\ker A=\s \{e_1\}$
and $Xe_1=0$ by the above lemma. Now we can use induction by writing 
$$
X=\begin{pmatrix} 

\; 0 & \vrule & \ast \;\\ 
\noalign{\hrule}  \\
\noalign{\vskip-13.7pt}
\; 0 & \vrule & X_1 \; \\
\end{pmatrix}
$$  
with $X_1\in M_{n-1}(K)$. It holds $X_1J(n-1)-J(n-1)X_1=X_1^p$ so that
$X_1$ is upper triangular by induction hypothesis. Hence $X$ is also
upper triangular.
\end{proof}

\begin{prop}
Let $p$ be an integer with $1<p<n$ and let $A=J(n)$, $X=(x_{i,j})\in M_n(K)$. 
Then $X$ is a matrix solution of $XA-AX=X^p$ if and only if
 
\begin{align}
x_{i,j} & = 0 \quad \text{for all} \quad 1\le j\le i\le n \label{j1} \\[0.2cm]
x_{i,j-1}-x_{i+1,j} & = \sum_{\ell_1=i+1}^{j-p+1} \sum_{\ell_2=\ell_1+1}^{j-p+2}
\cdots \sum_{\ell_{p-1}=\ell_{p-2}+1}^{j-1}x_{i,\ell_1}x_{\ell_1,\ell_2}
\cdots x_{\ell_{p-2},\ell_{p-1}}x_{\ell_{p-1},j}\label{j2}
\end{align}
for all $1\le i<j\le n$. 
\end{prop}

\begin{proof}
By proposition $\ref{prop12}$ we know that $X$ is upper triangular. 
Hence \eqref{j1} holds.
The equations \eqref{j2} follow by matrix multiplication. The $(i,j)$-th 
coefficient of $XA-AX$ is just the LHS of \eqref{j2} whereas the $(i,j)$-th
coefficient of $X^p$ is given by the RHS of \eqref{j2}. This may be seen
by induction.
\end{proof}

\begin{rem}
We can solve the polynomial equations given by \eqref{j2}
recursively. For $p\ge 3$ every $x_{i+1,j}$ can be expressed
as a polynomial in the free variables 
$x_{1,2},\ldots , x_{1,n}$ since the RHS of \eqref{j2} does not contain $x_{i+1,j}$.
For $p=2$ however it does contain $x_{i+1,j}$ for $\ell=i+1$.
In that case we rewrite the equations as follows.

\begin{align*}
x_{i+1,j}(1+x_{i,i+1}) - x_{i,j-1} + \sum_{\ell_1=i+2}^{j-1} x_{i,\ell}x_{\ell,j} & = 0 
\end{align*}
For $j=i+2$ we obtain $x_{i+1,i+2}(1+x_{i,i+1})=x_{i,i+1}$. This shows that
$1+x_{i,i+1}$ is always nonzero. It follows that every $x_{i+1,j}$ is a polynomial
in $x_{1,2},\ldots , x_{1,n}$ divided by a product of factors $1+kx_{1,2}$ also
being nonzero.
The formulas can be determined recursively. The first two are as follows:
\begin{align*}
x_{i+1,i+2} & = \frac{x_{1,2}}{1+ix_{1,2}}\\[0.3cm]
x_{i+1,i+3} & = \frac{x_{1,3}(1+x_{1,2})}{(1+ix_{1,2})(1+(i+1)x_{1,2})}\\
\end{align*}
\end{rem}

\begin{ex}
Let $n=5$ and $A=J(5)$. Then all matrix solutions
$X=(x_{ij})\in M_5(K)$ of $XA-AX=X^2$ are given by

$$
X=\begin{pmatrix}
       0 & x_{12} & x_{13} & x_{14} & x_{15} \\[0.2cm]
       0 & 0 & \frac{x_{12}}{1+x_{12}} & \frac{x_{13}}{1+2x_{12}} & 
\frac{(1+2x_{12})^2x_{14}-(1+x_{12})x_{13}^2}
{(1+x_{12})(1+2x_{12})(1+3x_{12})} \\[0.3cm]
       0 & 0 & 0 & \frac{x_{12}}{1+2x_{12}} & \frac{(1+x_{12})x_{13}}{(1+2x_{12})
(1+3x_{12})} \\[0.3cm]
       0 & 0 & 0 & 0 & \frac{x_{12}}{1+3x_{12}} \\[0.3cm]
       0 & 0 & 0 & 0 & 0
\end{pmatrix}.
$$
\end{ex}

\begin{cor}
Let $A=J(n)$. A special matrix solution of $XA-AX=X^2$ is given as follows:

$$
X_0:=\begin{pmatrix}
0 & \al & 0 & & & & \\
0 & 0 & \frac{\al}{1+\al} & 0 & & &\\
0 & 0 & 0 & \frac{\al}{1+2\al} & & &\\
0 & 0 & 0 & 0 &  & &\\
\vdots & \vdots & \vdots & \vdots & \ddots &  &\\
0 & 0 & 0 & 0 & \cdots & 0 & \frac{\al}{1+(n-2)\al} \\
0 & 0 & 0 & 0 & \cdots & \cdots & 0
\end{pmatrix}.
$$
\end{cor}

In some cases all matrix solutions of $XJ(n)-J(n)X=X^2$ are conjugated
to $X_0$.

\begin{prop}
Let $A=J(n)$ and $X=(x_{ij})\in M_n(K)$ be a matrix solution of 
$XA-AX=X^2$ with $x_{12}=\al\neq 0$. Then there exists an $S\in GL_n(K)\cap C(A)$
such that $X=SX_0S^{-1}$.
\end{prop}

\begin{proof}
We will prove the result by induction on $n$. The case $n=2$ is obvious.
For $A=J(n)$ we have $C(A)=\{c_1A^0+c_2A^1+\cdots + c_{n}A^{n-1}\mid c_i \in K\}$.
The matrix solution $X$ is strictly upper triangular by proposition $\ref{prop12}$.
We have
$$
X=\begin{pmatrix} 
\; X' & \vrule & \ast \;\\ 
\noalign{\hrule}  \\
\noalign{\vskip-13.7pt}
\; 0 \cdots 0 & \vrule & 0 \; \\
\end{pmatrix}
$$  
with $X'\in M_{n-1}(K)$. It is easy to see that $XA-AX=X^2$ implies
$X'A'-A'X'=X'^2$ with $A'=J(n-1)$. Hence by assumption there exists an
$S'\in GL_{n-1}(K)\cap C(A')$ such that

$$
S'X'S'^{-1}=X_0'
$$
where $X_0'$ is the special solution in dimension $n-1$.
We can extend $S'$ to a matrix $S_1\in GL_n(K)\cap C(A)$ as follows:

$$
S_1=
\begin{pmatrix} 
 & \vrule & s_n \\ 
\; S' & \vrule & \vdots \\
      & \vrule & s_2 \\
\noalign{\hrule}  \\
\noalign{\vskip-14pt}
\; 0 \cdots 0 & \vrule & s_1 \; \\
\end{pmatrix}
=
\begin{pmatrix}
s_1 & s_2 & \cdots & \cdots & \cdots & s_{n}\\
0 & s_1 & s_2 & \cdots & \cdots & s_{n-1}\\
\vdots & \vdots &  &  &  & \vdots \\
0 & 0 & 0 & \cdots & s_1 & s_2 \\
0 & 0 & 0 & \cdots & \cdots & s_1
\end{pmatrix}.
$$
One verifies that
\begin{align*}
S_1XS_1^{-1} & = 
\begin{pmatrix} 
 & \vrule & \ast \\ 
\; S' & \vrule & \vdots \\
      & \vrule & \ast \\
\noalign{\hrule}  \\
\noalign{\vskip-13.7pt}
\; 0 \cdots 0 & \vrule & s_1 \; 
\end{pmatrix}
\begin{pmatrix} 
 & \vrule & \ast \\ 
\; X' & \vrule & \vdots \\
      & \vrule & \ast \\
\noalign{\hrule}  \\
\noalign{\vskip-13.7pt}
\; 0 \cdots 0 & \vrule & 0 \; 
\end{pmatrix}
\begin{pmatrix} 
 & \vrule & \ast \\ 
\; S'^{-1} & \vrule & \vdots \\
      & \vrule & \ast \\
\noalign{\hrule}  \\
\noalign{\vskip-13.7pt}
\; 0 \cdots 0 & \vrule & s_1^{-1} \; 
\end{pmatrix}\\[0.5cm]
 & = 
\begin{pmatrix} 
 & \vrule & r_1 \\ 
\; X_0' & \vrule & \vdots \\
      & \vrule & r_{n-1} \\
\noalign{\hrule}  \\
\noalign{\vskip-13.7pt}
\; 0 \cdots 0 & \vrule & 0 \; 
\end{pmatrix}.
\end{align*}
The last matrix is not yet equal to $X_0$. It is however
a solution of $XA-AX=X^2$ by corollary $\ref{cor7}$. 
A short computation shows that this is true if and only if
\begin{align*}
0 & = r_i (1+(i-1)\al) \quad \text{for} \quad i=2,3,\ldots ,n-2 \\
 r_{n-1} & = \frac{\al}{1+(n-2)\al} 
\end{align*}
Hence we have $r_2 =\cdots = r_{n-2}=0$. 
It remains to achieve $r_1=0$. This is done by conjugating with
$$
S_2=
\begin{pmatrix}
1 & 0 & \cdots & r & 0\\
0 & 1 & \cdots & 0 & r\\
\vdots & \vdots &  & \vdots & \vdots \\
0 & 0 & \cdots & 1 & 0 \\
0 & 0 & \cdots & 0 & 1
\end{pmatrix}
$$
where 
$r=\frac{r_1(1+(n-2)\al)}{4\al^2}$.
Note that we have by assumption $\al\neq 0$. Let $S:=S_2S_1\in GL_n(K)\cap C(A)$.
We obtain $SXS^{-1}=X_0$.
\end{proof}

The solutions $X$ with $x_{12}=0$ need not be conjugated to $X_0$.

\begin{ex}
Let $A=J(8)$ and $X=\al A^4+\be A^5+\ga A^6+\de A^7$. Then 
$XA-AX=0=X^2$ and $SXS^{-1}=X$ for all $S\in GL_n(K)\cap C(A)$.
In particular $X$ is not conjugated to $X_0$.
\end{ex}

In general we may assume that $A$ is a Jordan block matrix with
Jordan blocks $J(r_1),\ldots ,J(r_k)$. If we have found solutions
$X_1,\ldots ,X_k$ to the equations $X_iJ(r_i)-J(r_i)X_i=X_i^p$ for
$i=1,2,\ldots ,k$ then $X=\diag(X_1,\ldots, X_k)$ is a solution
of $XA-AX=X^p$. However these are not the only solutions in general.
How can one determine the other solutions ?
One way would be to classify the matrix solutions of $XA-AX=X^p$ up to 
conjugation with the centralizer of $A$. 
First examples show that this classification will be complicated.
The following examples illustrate this for $A=\diag(J(2),J(2))$
and $p=2,3$.

\begin{ex}\label{gem1}
The matrix solutions of $XA-AX=X^2$ with
$$
A=\begin{pmatrix}
0 & 1 & 0 & 0 \\
0 & 0 & 0 & 0 \\
0 & 0 & 0 & 1 \\
0 & 0 & 0 & 0
\end{pmatrix}
$$
are given, up to conjugation with  $S\in GL_4(K)\cap C(A)$, by the
following matrices:

\begin{align*}
X_1 & =\begin{pmatrix}
0 & -1 & 0 & 0 \\
0 & 0 & 1 & 0 \\
0 & 0 & 0 & 1 \\
0 & 0 & 0 & 0
\end{pmatrix}, \quad
X_2 =\begin{pmatrix}
0 & -1 & 0 & 0 \\
0 & 0 & 0 & 1 \\
0 & 0 & 0 & 0 \\
0 & 0 & 0 & 0
\end{pmatrix}, \quad
X_{3,\al} =\begin{pmatrix}
0 & 0 & 1 & 0 \\
0 & 0 & 0 & \al \\
0 & 0 & 0 & 1-\al \\
0 & 0 & 0 & 0
\end{pmatrix}, \\[0.2cm] 
X_{4,\al,\be} & =
\begin{pmatrix}
0 & \al & 0 & 0 \\
0 & 0 & 0 & 0 \\
0 & 0 & 0 & \be \\
0 & 0 & 0 & 0
\end{pmatrix}, \quad
X_{5,\al}=
\begin{pmatrix}
0 & \al & 0 & 1 \\
0 & 0 & 0 & 0 \\
0 & 0 & 0 & \al \\
0 & 0 & 0 & 0
\end{pmatrix}.
\end{align*}
\end{ex}

The solutions $X_{4,\al,\be}=\diag (Y_1,Y_2)$ arise from the 
solutions of the equations $J(2)Y_i-Y_iJ(2)=X^2$. 
The solutions satisfying $X^2=0=XA-AX$ are given by
$X_{3,1}, X_{4,\al,\be}$ and $X_{5,\al}$.
The result is verified by an explicit computation. 
The centralizer of $A$
consists of matrices of the form
\begin{align*}
\begin{pmatrix}
s_1 & s_2 & s_5 & s_6 \\
0 & s_1 & 0 & s_5 \\
s_3 & s_4 & s_7 & s_8 \\
0 & s_3 & 0 & s_7
\end{pmatrix}
\end{align*}

whose determinant is given by $(s_1s_7-s_3s_5)^2$.

\begin{ex}
The matrix solutions of $XA-AX=X^3$ with $A=\diag(J(2),J(2))$
are given, up to conjugation with  $S\in GL_4(K)\cap C(A)$, by the
following matrices:

\begin{align*}
X_{1,\al,\be} & =\begin{pmatrix}
0 & 0 & \al & 0 \\
0 & 0 & 0 & \be \\
0 & \frac{\al-\be}{\al\be} & 0 & 0 \\
0 & 0 & 0 & 0
\end{pmatrix}, \quad
X_{2,\al} =\begin{pmatrix}
0 & \al & 1 & 0 \\
0 & 0 & 0 & 1 \\
0 & 0 & 0 & 0 \\
0 & 0 & 0 & 0
\end{pmatrix}, \\[0.2cm] 
X_{3,\al,\be} & =
\begin{pmatrix}
0 & \al & 0 & 0 \\
0 & 0 & 0 & 0 \\
0 & 0 & 0 & \be \\
0 & 0 & 0 & 0
\end{pmatrix}, \quad 
X_{4,\al} =
\begin{pmatrix}
0 & \al & 0 & 1 \\
0 & 0 & 0 & 0 \\
0 & 0 & 0 & \al \\
0 & 0 & 0 & 0
\end{pmatrix}.
\end{align*}
\end{ex}

The only solutions satisfying $X^3\neq 0$ are 
$X_{1,\al,\be}$ where $\al \neq \be $.

\section{The case $p=2$}

For $p=2$ our matrix equation is given by $X^2=XA-AX$. This equation
is a special case of the well known algebraic Riccati equation.
There is a large literature on this equation, see  \cite{LAR} and the references therein.
In particular, there is a well known result on the parametrization of solutions
of the Riccati equation using Jordan chains. The consequence is that matrix solutions
can be constructed by determining Jordan chains of certain matrices.
This does not mean, however, that we are able to solve the algebraic Riccati 
equation explicitly. The problem is only reformulated in terms of Jordan chains. 
Nevertheless this is an interesting approach. We will apply this result
to our special case and demonstrate it by an example. \\
The algebraic Riccati equation is the following quadratic
matrix equation \cite{LAR}

\begin{align*}
XBX + XA -DX -C & = 0
\end{align*}
where $A,B,C,D$ have sizes $n\times n, n\times m, m\times n$ and $m\times m$
respectively. Here $m\times n$ matrix solutions $X$ are to be found.
The special case $m=n$ and $B=-E,\, D=A,\, C=0$ yields $XA-AX-X^2=0$.

\begin{defi}
A {\it Jordan chain} of an $n\times n$ matrix $T$ is an ordered set
of vectors $x_1,\ldots x_r \in K^n$ such that $x_1\neq 0$ and for some
eigenvalue $\la$ of $T$ the equalities
\begin{align*}
(T-\la E)x_1 & = 0\\
(T-\la E)x_2 & = x_1\\
\vdots \hskip 1cm  & = \; \vdots \\
(T-\la E)x_r & = x_{r-1}
\end{align*}
hold. 
\end{defi}

The vectors $x_2,\ldots , x_r$ are called generalized eigenvectors
of $T$ associated with the eigenvalue $\la$ and the eigenvector $x_1$.
The number $r$ is called the length of the Jordan chain.

We call the $n$-dimensional subspace
\begin{align*}
G(X)& = \im  \begin{bmatrix} E \\ X \end{bmatrix} \subseteq K^{2n}
\end{align*}
the graph of $X$. Denote by $T\in M_{2n}(K)$ the matrix
\begin{align*}
T & = \begin{bmatrix} A & -E \\ 0 & A \end{bmatrix}.
\end{align*}
Then we have the following simple result \cite{LAR}:

\begin{prop}
For any $n\times n$ matrix $X$, the graph of $X$ is $T$-invariant
if and only if $X$ is a solution of $XA-AX=X^2$. 
\end{prop}

Representing the $T$-invariant subspace $G(X)$ as the linear span
of Jordan chains of $T$, we obtain the following result \cite{LAR}.

\begin{prop}
The matrix $X\in M_n(K)$ is a solution of $XA-AX=X^2$ if and only
if there is a set of vectors $v_1,\ldots, v_n \in K^{2n}$
consisting of sets of Jordan chains for $T$ such that
$$
v_i=\begin{bmatrix}  y_i \\ z_i
\end{bmatrix}
$$
where $y_i,z_i \in K^n$ and $(y_1, \ldots, y_n)$ forms a basis of $K^n$. 
Furthermore, if
\begin{align*}
Y & = \begin{bmatrix} y_1 & y_2 & \cdots & y_n \end{bmatrix} \in M_n(K),
\quad Z = \begin{bmatrix} z_1 & z_2 & \cdots & z_n \end{bmatrix} \in M_n(K),
\end{align*}
every matrix solution of $XA-AX=X^2$ has the form $X=ZY^{-1}$ for some set of
Jordan chains $v_1,\ldots , v_n$ for $T$, such that $Y$ is nonsingular.
\end{prop}

It follows that there is a one-to-one correspondence between the set of solutions
of $XA-AX=X^2$ and a certain subset of $n$-dimensional $T$-invariant
subspaces.

\begin{ex}
Let $A=\diag(J(2),J(2))\in M_4(K)$ and 
\begin{align*}
T & = \begin{bmatrix} A & -E \\ 0 & A \end{bmatrix} \in M_8(K)
\end{align*}
Then a set of Jordan chains for $T$ is given by
\begin{align*}
v_1 & = (2,0,0,0,0,0,0,0)^t \\
v_2 & = (0,1,0,0,-1,0,0,0)^t\\
v_3 & = (0,0,-1,0,0,-1,0,0)^t\\
v_4 & = (0,0,0,1,0,0,1,0)^t 
\end{align*}
We have $Tv_1=0,\, Tv_2=v_1,\, Tv_3=v_2$ and $Tv_4=v_3$. Then

\begin{align*}
X & = ZY^{-1} = \begin{pmatrix} 
0 & -1 & 0 & 0 \\
0 & 0 & -1 & 0 \\
0 & 0 & 0 & 1 \\
0 & 0 & 0 & 0
\end{pmatrix}
\begin{pmatrix} 
2 & 0 & 0 & 0 \\
0 & 1 & 0 & 0 \\
0 & 0 & -1 & 0 \\
0 & 0 & 0 & 1
\end{pmatrix}^{-1}=
\begin{pmatrix} 
0 & -1 & 0 & 0 \\
0 & 0 & 1 & 0 \\
0 & 0 & 0 & 1 \\
0 & 0 & 0 & 0
\end{pmatrix}
\end{align*}
is a matrix solution of $XA-AX=X^2$, see example $\ref{gem1}$.
\end{ex}

Note that Jordan chains of length three for $T$ are given by vectors 
of the form
\begin{align*}
v_1 & =(2a_1,0,2a_2,0,0,0,0,0)^t\\
v_2 & =(a_3,a_1,a_4,a_2,-a_1,0,-a_2,0)^t\\
v_3 & =(a_5,a_6,a_7,a_8,a_6-a_3,-a_1,a_8-a_4,-a_2)^t
\end{align*}

\section{Combinatorial formulas}

The matrix equation $XA-AX=X^p$ may be interpreted as a commutator rule
$[X,A]=X^p$. Successive commuting yields very interesting formulas
for $X^{\ell}A^m$ and $A^m X^{\ell}$, where $\ell,m\ge 1$. For $m=1$
the formulas are easy: we have $X^\ell A =AX^{\ell}+\ell X^{\ell+p-1}$ and
$AX^{\ell}= X^\ell A -\ell X^{\ell+p-1}$. For $m\ge 2$ these formulas become more
complicated. Finally we will prove a formula for $(AX)^{\ell}$. 
Although it is not needed for the study of solutions of our matrix equation, we would like to 
include this formula here. In fact, the commutator formulas presented here are important
for many topics in combinatorics. We are able to prove explicit 
formulas involving weighted Stirling numbers.

\begin{prop}
Let $p\ge 1$ and let $X,A\in M_n(K)$ satisfy the matrix equation $XA-AX=X^p$.
Then for $\ell,m \ge 1$ we have

\begin{align}\label{a1}
X^\ell A^m & =\sum_{k=0}^m a_{\ell}(k)\binom{m}{k}A^{m-k}X^{\ell+k(p-1)}
\end{align}
where $a_{\ell}(0)=1$ and $a_{\ell}(k)=a_{\ell, p}(k)=\prod_{j=0}^{k-1} [\ell+j(p-1)]$.
In particular we have

\begin{align}
X^\ell A & =AX^{\ell}+\ell X^{\ell+p-1}\\
X^\ell A^2 & = A^2X^{\ell}+2\ell AX^{\ell+p-1}+\ell(\ell+p-1)X^{\ell+2(p-1)}\\
\begin{split}
X^\ell A^3 & = A^3X^{\ell}+3\ell A^2X^{\ell+p-1}+3\ell (\ell+p-1)AX^{\ell+2(p-1)}
  +\ell (\ell+p-1) (\ell+2(p-1)) \\ 
& \quad \; X^{\ell+3(p-1)}.
\end{split}
\end{align}
For $p=2$ the formula simplifies to

\begin{align}\label{a3}
X^\ell A^m & =\sum_{k=0}^m k!\, \binom{m}{k} \binom{\ell+k-1}{\ell-1}
A^{m-k}X^{\ell+k}.
\end{align}
\end{prop}

\begin{proof}
For the case $m=1$  see \eqref{a2}. Now
\eqref{a1} follows by induction over $m$. Note that
$a_{\ell}(k+1)=a_{\ell}(k)(\ell+k(p-1))$.

\begin{align*}
X^{\ell}A^{m+1} & = \left(X^{\ell}A^m\right)A = \sum_{k=0}^m \binom{m}{k}
a_{\ell}(k)A^{m-k}\left(X^{\ell+k(p-1)}A\right)\\
 & = \sum_{k=0}^m \binom{m}{k}a_{\ell}(k)A^{m-k}\left(AX^{\ell+k(p-1)}+
(\ell+k(p-1))X^{\ell+(k+1)(p-1)}\right)\\
 & = \sum_{k=0}^m \binom{m}{k}a_{\ell}(k)A^{m+1-k}X^{\ell+k(p-1)}+
\sum_{k=0}^m \binom{m}{k}a_{\ell}(k)(\ell+k(p-1))A^{m-k}X^{\ell+(k+1)(p-1)}\\
 & =A^{m+1}X^{\ell}+ \sum_{k=1}^m \binom{m}{k}a_{\ell}(k)A^{m+1-k}X^{\ell+k(p-1)}\\
 & + \; \sum_{k=1}^m \binom{m}{k-1}\left[(\ell+(k-1)(p-1))a_{\ell}(k-1)\right]
A^{m+1-k}X^{\ell+k(p-1)}\\
 & + \; (\ell+m(p-1))a_{\ell}(m)X^{\ell+(m+1)(p-1)}\\
 & = \sum_{k=0}^{m+1}a_{\ell}(k)\binom{m+1}{k}A^{m+1-k}X^{\ell+k(p-1)}
\end{align*}

For $p=2$ we have $a_{\ell}(k)=\ell(\ell+1)\cdots (\ell+k-1)=k!\,
\binom{\ell+k-1}{k-1}$.
\end{proof}

In the same way one can prove the following result by induction:

\begin{prop}
Let $p\ge 1$ and let $X,A\in M_n(K)$ satisfy the matrix equation $XA-AX=X^p$.
Then for $\ell,m \ge 1$ we have

\begin{align}\label{a11}
A^mX^\ell & =\sum_{k=0}^m (-1)^k a_{\ell}(k)\binom{m}{k}X^{\ell+k(p-1)} A^{m-k}.
\end{align}

\end{prop}

\begin{prop}
Let $p\ge 2$ and let $X,A\in M_n(K)$ satisfy the matrix equation $XA-AX=X^p$.
Then we have for all $\ell\ge 1$

\begin{align}\label{a4}
(AX)^\ell & =\sum_{k=0}^{\ell-1}c(\ell,k)A^{\ell-k}X^{\ell+k(p-1)}
\end{align}
where the numbers $c(\ell,k)=c(\ell,k,p)$ are defined by the 
following recurrence relation for $1\le k\le \ell$.

\begin{align}
c(\ell,0) & =1 \\
c(\ell,\ell) & =0 \\
c(\ell+1,k) & =c(\ell,k)+\left[\ell+(p-1)(k-1)\right] c(\ell,k-1) \label{r1}
\end{align}
\end{prop}

\begin{proof}
We proceed by induction on $\ell$. Using \eqref{a2} we obtain

\begin{align*}
(AX)^{\ell+1} & = (AX)^{\ell}AX= \sum_{k=0}^{\ell-1} c(\ell,k)A^{\ell-k}
\left( X^{\ell+k(p-1)}A\right)X\\
 & = \sum_{k=0}^{\ell-1} c(\ell,k)A^{\ell-k}\left[ AX^{\ell+k(p-1)}+
(\ell+k(p-1))X^{\ell+(k+1)(p-1)}\right]X\\
 & = \sum_{k=0}^{\ell-1} c(\ell,k)A^{\ell+1-k}X^{\ell+1+k(p-1)}+
\sum_{k=0}^{\ell-1} c(\ell,k)(\ell+k(p-1))A^{\ell-k}X^{\ell+1+(k+1)(p-1)}.
\end{align*}

Using \eqref{r1} it follows
\begin{align*}
(AX)^{\ell+1} & =
 A^{\ell+1}X^{\ell+1}+ \sum_{k=1}^{\ell-1} c(\ell,k)A^{\ell+1-k}X^{\ell+1+k(p-1)}\\
 & + \; \sum_{k=1}^{\ell} c(\ell,k-1)(\ell+(k-1)(p-1))A^{\ell+1-k}X^{\ell+1+k(p-1)}\\
 & = A^{\ell+1}X^{\ell+1}+\sum_{k=1}^{\ell-1} c(\ell+1,k)A^{\ell+1-k}X^{\ell+1+k(p-1)}\\
 & + \, c(\ell,\ell-1)(\ell+(\ell-1)(p-1))AX^{\ell+1+\ell(p-1)}\\
 & = \sum_{k=0}^{\ell}c(\ell+1,k)A^{\ell+1-k}X^{\ell+1+k(p-1)}.
\end{align*}
\end{proof}

The integers $c(\ell,k,p)$ are uniquely determined.
The following table shows the values for 
$\ell=1,2,\ldots 6$ and $k=0,\ldots, \ell-1$

\vspace*{0.5cm}
\begin{Small}
\begin{center}
\begin{tabular}{c|c|c|c|c|c|c}
 $\ell \setminus k$ & 0 & 1 & 2 & 3 & 4 & 5 \\
\hline
$1$ & $1$ &  &  &  &  &  \\
$2$ & $1$ & $1$ & & & & \\
$3$ & $1$ & $3$ & $p+1$ & & & \\
$4$ & $1$ & $6$ & $4p+7$ & $(p+1)(2p+1)$ & & \\
$5$ & $1$ & $10$ & $5(2p+5)$ & $5(p+1)(2p+3)$ & $(p+1)(2p+1)(3p+1)$ & \\
$6$ & $1$ & $15$ & $5(4p+13)$ & $15(p+2)(2p+3)$ & $(p+1)(36p^2+70p+31)$ &
$(p+1)(2p+1)(3p+1)(4p+1)$\\
\end{tabular}
\end{center}
\end{Small}
\vspace*{0.5cm}

It is possible to find an explicit formula for the $c(\ell,k,p)$.

\begin{prop}
For $\ell\ge 1$ and $0\le k\le \ell-1$ we have

\begin{align}\label{c1}
c(\ell,k,p) & = (p-1)^{k-\ell +1}\sum_{r=1}^{\ell-k}\frac{(-1)^{r-1}}{(r-1)!\,
(\ell-k-r)!}\; \prod_{j=1}^{\ell -1}[pj+(1-p)r].
\end{align}
For $p=2$ the formula reduces to

\begin{align}
c(\ell,k,2) & = \binom{\ell+k-1}{2k} \prod_{j=1}^{k}(2j-1).
\end{align}

\end{prop}

\begin{proof}
Let $S(n,k)=S(n,k, \la | \theta)$ denote the  weighted degenerated Stirling numbers
for $n,k \ge 1$, see  \cite{HOW}.
They are given by
\begin{align}
S(n,n) & = 1 \label{s1}\\
S(n,0) & = \prod_{j=0}^{n-1}(\la -j\theta)\label{s2}\\
S(n+1,k) & = (k+\la -\theta n)S(n,k)+S(n,k-1). \label{s3}
\end{align}
\vspace*{0.3cm}
They satisfy the following explicit formula (see $(4.2)$ of \cite{HOW}):

\begin{align}\label{s4}
S(n,k,\la | \theta) & = \sum_{r=0}^k \frac{(-1)^{k+r}}{r!\, (k-r)!}\, \prod_{j=0}^{n-1}
(\la +r-j\theta).
\end{align}
We can rewrite the recurrence relation \eqref{r1} for the numbers $c(\ell,k)$
as follows. If we substitute $k$ by $\ell-k+1$ then we obtain

\begin{align*}
c(\ell+1,\ell-k+1) & = c(\ell,\ell-k+1)+ \left[\ell+(p-1)(\ell-k)\right] c(\ell,\ell-k).
\end{align*}

Here $\ell-k+1$ runs through $1,2,\ldots \ell$ if $k$ does.
Now set $c(\ell,\ell-k)=s(\ell,k)(1-p)^{\ell-k}$. Then the above recurrence
relation implies that

\begin{align*}
s(\ell,\ell) & = c(\ell,0)=1\\
s(\ell,0) & = c(\ell,\ell)(1-p)^{-\ell}=0=\prod_{j=0}^{\ell-1}pj\\
s(\ell+1,k)& =s(\ell,k-1)+\left[ k+\left( \frac{p}{1-p}\right) \right]s(\ell,k).
\end{align*}
\vspace*{0.3cm}
Comparing this with \eqref{s1},\eqref{s2},\eqref{s3} we see that 
$s(\ell,k)=S(\ell,k, \la| \theta)$ for $\la=0$ and $\theta= \frac{p}{p-1}$.
So they are indeed degenerated Stirling numbers. Applying the formula \eqref{s4} 
to $c(\ell,k)=s(\ell,\ell-k)(1-p)^{k}$ we obtain
\begin{align*}
c(\ell,k,p) & = (p-1)^k\sum_{r=1}^{\ell-k}\frac{(-1)^{r-1}}{(r-1)!\,
(\ell-k-r)!}\; \prod_{j=1}^{\ell -1}\left[\frac{pj}{p-1}-r\right].
\end{align*}
This shows \eqref{c1}. For $p=2$ one can obtain
a much easier formula. It is however easier to derive this formula
not from \eqref{c1} but rather by direct verification of the recurrence
relation.
\end{proof}

\begin{rem}
We have the following special cases:

\begin{align*}
c(\ell,1,p) & = \frac{\ell(\ell-1)}{2}\\[0.3cm]
c(\ell,2,p) & = \frac{\ell(\ell-1)(\ell-2)(3\ell+4p-5)}{24}\\[0.3cm]
c(\ell,\ell-1,p) & = (1+p)(1+2p)\cdots (1+(\ell-2)p).
\end{align*}
\end{rem}

\end{document}